\documentclass[a4paper,11pt]{article}
\usepackage{array}
\usepackage{theorem}
\usepackage{amsmath,amscd,amssymb} 
\usepackage{latexsym}
\usepackage{epsfig}
\usepackage[all]{xy}
\theorembodyfont{\sl}

\newtheorem{lemma}{Lemma}[section]

\newtheorem{proposition}[lemma]{Proposition}
\newtheorem{theorem}[lemma]{Theorem}
\newtheorem{corollary}[lemma]{Corollary}

\newtheorem{definition}[lemma]{Definition}

\renewcommand{\AA}{\mathbb A}

\newcommand{\FF}{\mathbb F}

\newcommand{\KK}{\mathbb K}

\newcommand{\PP}{\mathbb P}

\newcommand{\cA}{\mathcal A}

\newcommand{\cC}{\mathcal C}

\newcommand{\cE}{\mathcal E}
\newcommand{\cF}{\mathcal F}

\newcommand{\cI}{\mathcal I}

\newcommand{\cL}{\mathcal L}
\newcommand{\cM}{\mathcal M}

\newcommand{\cO}{\mathcal O}

\newcommand{\cQ}{\mathcal Q}
\newcommand{\cR}{\mathcal R}

\newcommand{\cX}{\mathcal X}

\newcommand{\tensor}{\otimes}

\renewcommand{\Tilde}{\widetilde}

\renewcommand{\Bar}{\overline}

\newcommand{\cross}{\times}

\newcommand{\imic}{\cong}

\newcommand{\GL}{\mathop{\mathrm {GL}}\nolimits}

\newcommand{\PGL}{\mathop{\mathrm {PGL}}\nolimits}

\newcommand{\Proj}{\mathop{\null\mathrm {Proj}}\nolimits}

\newcommand{\Gr}{\mathop{\mathrm {Gr}}\nolimits}
\newcommand{\Jac}{\mathop{\mathrm {Jac}}\nolimits}
\newcommand{\Pic}{\mathop{\mathrm {Pic}}\nolimits}
\newcommand{\Sing}{\mathop{\mathrm {Sing}}\nolimits}
\newcommand{\Spec}{\mathop{\mathrm {Spec}}\nolimits}

\newcommand{\Hilb}{\mathop{\mathrm {Hilb}}\nolimits}

\newcommand{\coker}{\mathop{\mathrm {coker}}\nolimits}

\newcommand{\bl}{\mathop{\mathrm {Bl}}\nolimits}

\newcommand{\bdm}{\mathbf m}

\newcommand{\bdr}{\mathbf r}

\newcommand{\gothm}{\mathfrak m}
\newcommand{\gothp}{\mathfrak p}
\newcommand{\gothq}{\mathfrak q}

\newcommand{\qx}[1]{\phantom{}^t\underline{x}{#1}\underline{x}}
\newcommand{\qedsymbol}{\mbox{$\Box$}}
\newcommand{\qed}{\unskip\nobreak\hfil\penalty50\hskip1em\hbox{}\nobreak
\hfill\qedsymbol\parfillskip=0pt\finalhyphendemerits=0}
\newenvironment{proof}{\begin{ProofwCaption}{Proof}}{\end{ProofwCaption}}
\newenvironment{ProofwCaption}[1]
 {\addvspace\theorempreskipamount \noindent{\it #1.}\rm}
 {\qed \par \addvspace\theorempostskipamount}

\begin{document}

\title{The moduli space of \'etale double covers of genus 5 curves is
  unirational} 
\author{E.~Izadi, M.~Lo~Giudice and G.K.~Sankaran}
\date{}
\maketitle

\begin{abstract}
  We show that the coarse moduli space $\cR_5$ of \'etale double
  covers of curves of genus~$5$ over the complex numbers is
  unirational.  We give two slightly different arguments, one purely
  geometric and the other more computational.
\end{abstract}

\section{Introduction}\label{intro}

The coarse moduli space $\cR_g$ of \'etale double covers of genus $g$
curves is sometimes referred to as the Prym moduli space. It can be
thought of as the moduli space of curves $C\in\cM_{g}$ equipped with a
nontrivial line bundle $\cL$ whose square is trivial. Thus $\cR_g$ is
also equipped with a morphism to $\cM_g$, which is a finite cover of
degree $2^{2g}-1$.  It has been extensively studied for small values
of $g$. In particular Donagi showed in \cite{Donagi} that $\cR_6$ is
unirational. Other proofs of the unirationality of $\cR_6$ were given
by by Verra \cite{Verra1} and by Mori and Mukai \cite{MoriMukai}. 

For $g\le 6$ the Prym map $p_g\colon \cR_g \to \cA_{g-1}$, which
associates to an \'etale double cover $\tau\colon \tilde C \to C$ the
Prym variety $P(\tau)=\coker(\tau^*\colon \Jac C \to \Jac \tilde C)$,
is dominant.  It therefore follows from Donagi's result that the
moduli space $\cA_5$ of principally polarised abelian $5$-folds is
unirational. 

Catanese \cite{Catanese} showed that $\cR_4$ is rational. Moreover
$\cR_1=X_0(2)$ is rational, and $\cR_3$ is birational to the moduli
space of bielliptic curves of genus $4$ which is known to be rational
(\cite{BardelliDelCentina}).  Also, $\cR_2$ is dominated by $\cA_2(2)$
which is rational (it is modelled by the Segre cubic). Hence $\cR_2$
is at least unirational. 

Clemens \cite{Clemens1} showed that $\cA_4$ is unirational, but using
intermediate Jacobians, not Prym varieties. In the introduction to
\cite{Catanese} it is stated that $\cR_5$ is unirational and a
reference is given to \cite{Clemens1}, then unpublished; but as far as
we can determine no proof is given there or anywhere else. 

In this paper we fill this gap by proving (Theorem~\ref{R5.unirat})
that $\cR_5$ is indeed unirational. We work over an algebraically
closed field $\KK$ of characteristic different from~$2$ (except in
Section~\ref{cancurves}). 

The basic construction used in our proof is to be found
in~\cite{Clemens1}.  If $X$ is a quartic surface in $\PP^3$ with
six ordinary double points at $P_0,\dots,P_5$ and no other
singularities we define $C_X$ to be the discriminant of the
projection from $P_0$. Generically it is a $5$-nodal plane sextic (hence of
genus~$5$), with an everywhere tangent conic coming from the tangent
cone to $X$ at $P_0$.  The quartic double solid branched along~$X$ has
(after blowing up) the structure of a conic bundle over $\PP^2$ with
discriminant curve $C_X$. Blowing up in the remaining five points
yields a conic bundle over a degree $4$ del Pezzo surface, and the
discriminant is the canonical model $\Tilde C_X$ of $C_X$.  This
determines a connected \'etale double cover of $\Tilde C_X$. 

The space $\cQ$ of quartic surfaces in $\PP^3$ with
six isolated ordinary double points, one of which is marked, is
unirational. This is well-known and quite easy to prove: see
Corollary~\ref{Q.unirat}. 

The construction above defines a morphism from the
unirational variety $\cQ$ to $\cR_5$, which is in turn endowed with a
finite (in fact $1024$-to-$1$) natural projection to $\cM_5$. Since
$\cR_5$ is irreducible, to prove the unirationality it is now enough
to prove that the map to $\cM_5$ is generically surjective. 

We present two different proofs that the map $\theta\colon \cQ \to \cM_5$
is dominant. One method exploits the special geometry of the family, 
using ideas of Donagi as worked out in \cite{Izadi}. We show by a dimension count
that the general genus~$5$ curve does have a plane model as a $5$-nodal sextic with 
an everwhere tangent conic, and then show how to recover the quartic in $\PP^3$ as 
a certain image of the double cover of $\PP^2$ branched along the sextic. 

The other approach, which was used in \cite{LoGiudice}, is
computational, and is applicable to any family of $5$-nodal
sextics. It uses the fact that $\Tilde C_X$ is a canonical curve, and
reduces the question of surjectivity of the Kodaira-Spencer map to
computing the rank of a certain matrix. This can then be verified at a
test point. 

Acknowledgements: Part of this paper forms part of the Ph.D. 
thesis~\cite{LoGiudice} of MLG, supervised by GKS and supported by
EPSRC (UK) and Universit\`a di Milano. We also wish to thank Sandro
Verra for telling us about the rationality of $\cR_3$. After we had
completed this paper we learned
about the paper \cite{Verra2} by Verra where he proves the
unirationality of the universal abelian variety over $\cA_4$. The
unirationality of $\cR_5$ also follows from his result. Note that our
proofs and his are entirely different. 

\section{Nodal quartics and nodal curves}\label{quartics}

The equation of a quartic surface $X\subset\PP^3$ with an isolated
ordinary double point at $P_0 = (0:0:0:1)$ is $F=u_2x_3^2 + 2u_3x_3 +
u_4=0$, where $u_d$ is a form of degree $d$ in $\KK[x_0,x_1,x_2]$ and
the quadratic form $u_2$ is non-degenerate.  The projection $\pi\colon
X\setminus{P_0} \to \PP^2$ from $P_0$ is induced by the homomorphism
\begin{equation*}
r\colon \KK[x_0,x_1,x_2] \to \KK[x_0,\dots, x_3]/(F),
\end{equation*}
sending $x_i$ to $x_i+(F)$ for $i=0,1,2$. If $X$ is general, then any
line in $\PP^3$ through $P_0$ intersects $X$ in at most two other
points so $\pi$ is a quasi-finite morphism. It is finite away from the
image $(u_2= 0) \subseteq \PP^2$ of the tangent cone to $X$ at $P_0$. 

We define the plane curve $C_X$ to be the locus of lines through $P_0$
tangent to $X$ away from $P_0$.  The following is easy to prove: see
\cite[Lemma~5.1]{Kreussler1}. 

\begin{proposition}\label{discrim_sextic}
$C_X$ is a plane curve of degree six given by $u_3^2-u_2u_4=0$. 
Furthermore if $Q \in X$ is a singular point (different from $P$),
then $\pi(Q)$ is a singular point of $C_X$. 
\end{proposition}

The locus of points in $\PP^2$ whose reduced fibre under $\pi$
consists of only one point is not irreducible. There are two
components, the curve $C_X$ and the conic $u_2=0$. 

\begin{definition}
A conic $V\subset \PP^2$ is called a contact conic of $C_X$ if $V$
cuts on $X$ a divisor which is divisible by $2$. 
\end{definition}

\begin{corollary}\label{contact_conic}
  Let $\bl_{P_0}(X) \to X$ be the blow-up of the surface $X$ at the
  point $P_0$. Then the unique morphism $\bl_{P_0}(X) \to \PP^2$ that
  commutes with $\pi$ is a double cover of $\PP^2$ branched along $C_X$. The
  image of the exceptional divisor in $\bl_{P_0}(X)$ is the contact
  conic of $C_X$ defined by the equation $u_2=0$. 
\end{corollary}
\begin{proof}
  It is essentially enough to observe that the tangent cone of $X$ at
  $P_0$ is defined by the equation $u_2 = 0$ in $\PP^3$.  The
  exceptional curve inside $\bl_{P_0}(X)$ corresponds to the set of
  lines in the tangent cone. To see that the conic $u_2=0$ is a
  contact conic of $C_X$ simply look at the ideal $(u_2,u_3^2-u_2u_4)
  = (u_2,u_3^2)$ in $\KK[x_0,x_1,x_2]$.  In particular this means that
  the points of contact are given by $u_2=u_3=0$
\end{proof}
Next we describe the singular locus of $C_X$. 
\begin{lemma}\label{sings_CX}
  Let $Y_d\subset\PP^3$ be the cone of vertex $P_0$ defined by the form
  $u_d$, and let $Q \in X$ be any point different from $P_0$ such that
  $\pi(Q) \in C_X$. Then $\pi(Q)\in\Sing C_X$ if and only if
  $Q\in(\Sing X)\cup (Y_2 \cap Y_3 \cap Y_4)$. 
\end{lemma}
\begin{proof}
If $\gothq$ is the homogeneous ideal of $Q$, then $Q \in X$ means $F \in
\gothq$ and $\pi(Q) \in C_X$ means $u_3^2-u_2u_4 \in \gothq$. Then,
from the equality
\begin{equation*}
u_2F = (u_2x_3+u_3)^2-(u_3^2-u_2u_4)
\end{equation*}
we obtain immediately $u_2x_3+u_3\in\gothq$. Now
$\pi(Q)$ is a singular point if and only if $u_3^2-u_2u_4 \in
\gothq^2$, and this happens if and only if $u_2F\in\gothq^2$, which
means that either $u_2\in\gothq$ or $F \in \gothq^2$. In the
latter case $Q$ is a singular point of $X$; in the former we also have
$u_3\in\gothq$, since $u_2x_3+u_3 \in\gothq$, and $u_4\in\gothq$ since
$2u_3x_3+u_4\in\gothq$. 
\end{proof}
\begin{proposition}\label{general_CX}
  For a general quartic surface $X$ with an isolated double point $P_0$, 
  the singularities of $C_X$ all come from
  singularities of $X$. 
\end{proposition}

\begin{proof}
  Suppose that $P_0\neq Q\in Y_2 \cap Y_3 \cap Y_4$. Then $\pi(Q)\in
  (u_2=u_3=u_4=0)$, which is empty for general $(u_2,u_3,u_4)$. 
\end{proof}

If $X$ is a quartic surface with at least one ordinary double point $P_0$,
we let $p\colon \Lambda_X \to\PP^3$ be the double cover branched along $X$
and let $W_X=\bl_{P_0}(\Lambda_X)$ be the blow-up of $\Lambda_X$ at $P_0$. 

\begin{proposition}\label{conic_bundle}
Let $X$, $P_0$, $\Lambda_X$ and $W_X$ be as above. The unique
morphism $f$ that makes the diagram
\begin{equation*}
\xymatrix{{W_X} \ar[r] \ar[drr]_{f} & {\Lambda_X} \ar[r]^{p}
         & {\PP^3} \ar@{-->}[d]^{\pi} \\ & & {\PP^2}   }
\end{equation*}
commute is a conic bundle over $\PP^2$, and the curve $C_X$ is the
locus of points whose fibre is a degenerate conic. 
\end{proposition}
\begin{proof} 
This is shown in \cite[p.222]{Clemens1}: there is a
more detailed version in \cite[Section 2]{Kreussler2}. 
\end{proof}
Kreussler \cite[Section 2]{Kreussler2} gives an explicit equation for
$W_X$ as a divisor inside a $\PP^2$ bundle over $\PP^2$. Put
$\cE=\cO\oplus\cO(-1)\oplus\cO(-2)$ over $\mathbb{P}^2$, and consider
$p\colon\PP(\cE)\to \PP^2$. Here $\PP (E)$ is the projective bundle of
hyperplanes in the fibres of $E$. Let $z_k \in H^0(\PP^2,\cE(k))$,
$k=0$, $1$, $2$ be {\em constant} non-zero global sections and define
the divisor $W_X\subset\PP(\cE)$ by
\begin{equation*}
-z_2^2+z_1^2u_2+2z_1z_0u_3+z_0^2u_4=0
\end{equation*}
(the left-hand side is a section of $\cO_{\PP(\cE)}(2)\otimes
p^*\cO_{\PP^2}(4)$). 

\begin{lemma}\label{canonical_curve}
 Let $C$ be a plane sextic curve whose only singularities are five
nodes in linear general position. Let $\sigma_C\colon \Tilde{\PP^2} \to
\PP^2$ denote the blow-up of $\PP^2$ in these five points. Then
$\Tilde\PP^2$ is a degree $4$ Del Pezzo surface and the anticanonical
embedding of $\Tilde\PP^2$ in $\PP^4$ realises the strict transform
$\Tilde C$ of $C$ as a smooth canonically embedded curve of genus~$5$. 
\end{lemma}
\begin{proof} 
  This is well known. That $\Tilde C$ is canonically embedded follows
  from a simple adjunction computation: if $E$ is the exceptional
  divisor and $H$ is the class of a line in $\PP^2$, then in
  $\Pic\Tilde \PP^2$ we have $\Tilde
  C=\sigma_C^*(6H)-2E=-2K_{\Tilde\PP^2}$. 
\end{proof}

Next we consider a more special case. Again $X\subset\PP^3$ is a
quartic surface but now $X$ is to have six isolated ordinary double
points, $P_0,\dots,P_5$ and no other singularities. With respect to
$P_0$ we will also assume $X$ to be general, in the precise sense of
Proposition~\ref{general_CX}. Under these hypotheses
$C_X$ is a plane sextic with precisely five nodes, at $\Bar
P_i=\pi(P_i)$, $1\le i\le 5$. We also assume that the $\Bar P_i$
are in linear general position. 

Let $f\colon W_X\to \PP^2$ be the conic bundle as in
Proposition~\ref{conic_bundle}, and let $\sigma_{C_X}\colon
\Tilde\PP^2\to \PP^2$ and $\Tilde C_X$ be as in
Lemma~\ref{canonical_curve}. Let $\Sigma$ be the surface
$f^{-1}(C_X)\subset W_X$ and put $S=\Sigma\times_{C_X}\Tilde C_X$,
with $\tilde f\colon S\to \Tilde C_X$ the projection. 

Let $\nu\colon \Tilde{S} \to S$ be the normalisation of $S$, and
consider the Stein factorisation of $\tilde f\circ\nu \colon \Tilde
S\to \Tilde{C}_X$
\begin{equation*}
\xymatrix{ {\Tilde{S}} \ar[r]^{\nu} \ar[d]_{f'} & {S} \ar[d]^{\Tilde{f}} \\
           {\Gamma} \ar[r]^{g} & {\Tilde{C}_X} } 
\end{equation*}
so $f'$ is a projective morphism with connected fibres and $g$ is a
finite morphism. 

\begin{proposition}\label{double_cover}
  For general $X$, the finite morphism $g\colon \Gamma\to \Tilde C_X$
  is an \'etale degree~$2$ map between smooth connected curves, naturally
  associated with the pair $(X,P_0)$. 
\end{proposition}
\begin{proof}
    The cover $g$ is unbranched of degree $2$ because the
    restriction of the conic bundle $W_X$ to the curve $C_X$ consists
    of a fibration by pairs of distinct lines. To see this recall the
    equation for the conic bundle $W_X$. The preimage of
    $x=(x_1:x_2:x_3)\in \PP^2$ is given by
\begin{equation*}
-z_2^2+z_1^2u_2(x)+2z_1 z_0u_3(x)+z_0^2u_4(x)=0,
\end{equation*}
and this has rank $2$ since $u_2$, $u_3$ and
$u_4$ never vanish simultaneously. 

It remains to check that $g$ is nontrivial, that is, that $\Gamma$ is
connected.  In characteristic zero this follows from the fact that the
Prym variety $P(g)$ is isomorphic to the intermediate Jacobian of the
conic bundle~\cite[Section 2]{Beauville3}. If the double cover were
trivial, then $P(g)$ would have dimension~$5$, which is impossible. To
extend this to the case of characteristic $p\neq 2$, we observe that
the quartic equation lifts to characteristic zero and in that case the
double cover is nontrivial, as we have just seen. Therefore the double
cover in positive characteristic is connected.
\end{proof}

\section{Moduli of curves}\label{moduli}

The functor $\bdr_g$ ($g\ge 2$) given by families of smooth projective
curves of genus $g$ with a connected \'etale double cover is coarsely
represented by an irreducible quasi-projective scheme $\cR_g$: see
\cite[\S{6}]{Beauville2}. The dimension of this moduli space is
$3g-3$, the same as the dimension of $\cM_g$, the extra data being one
of the $2^{2g}-1$ nontrivial $2$-torsion points in the Jacobian of
$C$. So forgetting this defines a natural transformation between
$\bdr_g$ and $\bdm_g$, and thus a morphism $\cR_g \to
\cM_g$ with finite fibres. 

Fix five points $P_1,\dots,P_5$ in $\PP^3$ in linear general
position. 
Let $\cQ$ be the space of quartic surfaces in $\PP^3$ with
ordinary double points at the $P_i$ and one additional ordinary double point distinct from the $P_i$. 

\begin{proposition}\label{fine_moduli}
$\cQ$ is an irreducible
 locally closed subscheme of the Hilbert scheme of
quartic surfaces in $\PP^3$ hence inherits a universal family of quartics from the Hilbert scheme. 
\end{proposition}

\begin{proof}
The Hilbert scheme of quartic surfaces in
$\PP^3$ is $\Hilb_{2m^2+2}(\PP^3)=\PP\big(H^0(\PP^3,\cO(4))\big)$. 

Let $\gothp_i$ be the homogeneous ideal of $P_i$ in $\PP^3$. The set of quartic
surfaces in $\PP^3$ with five double points at the $P_i$ is the closed subscheme of the
Hilbert scheme given by $\PP(I_4)$ where $I =
\bigcap_{i=0}^4\gothp_i^2$ is the ideal of the five double points. In
other words it is $\PP H^0(\PP^3,\cI(4))$, where we denote $\cI$ the
sheaf of ideals defined by $I$. Note that $h^0(\PP^3,\cI(4))=15$. 

Going on, we ask now for a sixth double point. We take the product $\PP(I_4)
\times \PP^3$ and consider the closed subscheme $B_0$ defined as
\begin{equation*}
B_0 :=\big\{(F,\gothp_0) \,\,\big|\,\, F\in\gothp_0^2 \big\} =
 \left\{ (F,P_0) \,\,\left|\,\, F(P_0) =
\frac{\partial{F}}{\partial{x_i}}(P_0) = 0
 \right\}\right. . 
\end{equation*}
The projection
\[
B_0\longrightarrow \PP^3
\]
is surjective and the fibres are linear spaces. Hence $B_0$ is
irreducible and rational.  Projecting onto $\PP(I_4)$, we consider the
scheme theoretic image of $B_0$, which is a closed subscheme $B$ of
$\PP(I_4)$, and observe that in the universal factorisation
\begin{equation*}
\xymatrix{ {B_0} \ar[dr] \ar[r]^{\beta} & {B} \ar[d] \\ & {\PP(I_4)} }
\end{equation*}
the dominant morphism $\beta$ is also proper, because it is the
external morphism of a composition which is proper. As a consequence
$\beta$ is surjective and the scheme theoretic image of $B_0$
coincides with the set theoretic image. Observe also that $B$ is
irreducible. 
However $B$ contains all the possible degenerations of a quartic
surface with six double points, while we are interested in those
surfaces with ordinary double points and no other singularities. 
But this is clearly an open condition, so we have proved that $\cQ$ is
an open subset of an irreducible closed subset of the Hilbert scheme. 
\end{proof}

\begin{corollary}\label{Q.unirat}
$\cQ$ is unirational of dimension~$13$. 
\end{corollary}

\begin{proof}
  We have seen in the proof of Proposition~\ref{fine_moduli}
  above that there is a dominant rational map from a closed
  subscheme $B_0$ of the product $\PP(I_4) \times \PP^3$ to $\cQ$,
  given by the projection onto $\PP(I_4)$. Now focus on the other
  projection and observe that a generic point $P_0\in\PP^3$ defines
  four independent conditions on the linear space $H^0(\PP^3,\cI(4))$. 
  Over an open subset $U\subset\PP^3$ these conditions define a
  vector bundle $E$ of rank $11$. To see that $U$ is not empty it is
  enough to fix a sixth point in $\PP^3$ and compute the Hilbert
  function of the ideal $J = \bigcap_{i=0}^5\gothp_i^2$, which we may
  easily do with {\it Macaulay}. The projective space bundle $B_1$
  over $U$ associated to $E$ is a rational variety embedded in $B_0$
  as a dense open subset, and mapping dominantly onto $\cQ$, which is
  therefore at least unirational. 

One may check by computing the differential of the projection at one
point that the dimension of $\cQ$ is $13$ (it is irreducible because
it is an open subset of an irreducible variety). 
\end{proof}

\begin{proposition}\label{classifying_map}
There exists a morphism of schemes $\varrho\colon \cQ \to \cR_5$ given
by the constructions above that associates to any nodal quartic
surface in $\PP^3$ a nodal sextic plane curve with a double cover. 
\end{proposition}

\begin{proof}
  We must globalise our earlier constructions. This is a standard
  gluing argument. Suppose first that the base scheme $B=\Spec A$ for
some ring $A$. Associate to the scheme $\cX = \Proj
  A[x_0,\dots,x_3]/(u_2x_3^2+2u_3x_3+u_4)$ the plane curve over $A$
  defined by the equation $u_3^2-u_2u_4$. This association is
  natural, in that it commutes with pull-backs. Indeed for any
  homomorphism of rings $A \to A'$ the pull-back of $\cX$ is given by $\Proj
(A\otimes{A'})[x_0,x_1,x_2]/(u_3^2-u_2u_4)$, and this is the same graded ring 
one would obtain by first pulling back
the family of surfaces and then applying the correspondence. 
\end{proof}

We want to use the morphism $\varrho$ to prove the
unirationality of $\cR_5$. To do so we must show that $\varrho$ is
dominant. We can simplify the problem by taking advantage of the
irreducibility of $\cR_g$ and $\cM_g$. We have a commutative diagram
\begin{equation*}
\xymatrix{ {\cQ} \ar[r]^{\varrho} \ar[dr]_{\theta} &
  {\cR_5} \ar[d]^{\eta} \\ & {\cM_5}  }
\end{equation*}
where $\eta$ forgets the double cover. 

\begin{lemma}\label{moduli_map}
The morphism $\varrho$ is dominant if and only if $\theta$ is dominant. 
\end{lemma}

\begin{proof}
  Since $\eta$ is a dominant morphism between irreducible spaces, it
  is immediate that $\theta$ is dominant if $\varrho$ is. Conversely,
  if $\varrho$ is not dominant then the scheme theoretic image
  $\varrho(\cQ)\subset\cR_5$ has dimension less than $\dim{\cR_5}
  = 12$, because $\cR_5$ is irreducible, so $\dim\theta(\cQ)<12$
  also, so $\theta$ is not dominant. 
\end{proof}

\section{Reconstructing the double solid}\label{doublesolid}

In this section we give a proof that $\theta\colon \cQ\to\cM_5$ is
dominant by making use of the special geometry of the family $\cQ$. We
show how to reconstruct the quartic double solid from a suitable plane
sextic model of a sufficiently general genus $5$ curve. 

\begin{lemma}\label{good_sextic}
  For a general $C\in \cM_5$ there exists a birational map $C\to\Bar
  C\subset \PP^2$ to a plane $5$-nodal sextic $\Bar C$, such that
  $\Bar C$ admits a contact conic $V\subset\PP^2$ meeting $\Bar C$
  only at six distinct smooth points of $\Bar C$. Furthermore, a
  general $C\in\cM_5$ has a one-parameter family of such birational
  plane models. 
\end{lemma}
\begin{proof}
  By the Kempf-Kleiman-Laksov theorem (\cite[V (1.1)]{ACGH}), a
  general $C\in \cM_5$ has a $2$-dimensional family $G_6^2(C)$ of $g_6^2$s and
  hence of birational models $\Bar C$ as a plane sextic. The general
  such sextic, for any given general $C$, has five nodes. 

For fixed $C$, the image of the map $|\cO_{\PP^2}(2)|\times
G_6^2(C)\imic \PP^5\times\PP^2 \to \Hilb_{12}\PP^1$ given by $(V,\Bar
C)\mapsto V\cap\Bar C$ intersects the codimension $6$ locus consisting
of subschemes with multiplicity at least $2$ at each point in a
variety of dimension $1$. In particular this intersection is
nonempty. If we take $\Bar C$ defined by projection from a node of a
general $6$-nodal quartic surface in $\PP^3$, i.e. we take $F$ as in
Section~\ref{quartics} with $u_d$ general, and $V=(u_2=0)$, we obtain
a pair $(V,\Bar C)$ whose intersection is six distinct smooth points. 

Thus for a general genus $5$ curve $C$, there is a $1$-dimensional
family of plane $5$-nodal sextic models $\Bar C$ of $C$, each having a
contact conic meeting $\Bar C$ at six distinct smooth points of $\Bar
C$. 
\end{proof}
\begin{proposition}\label{rebuild_double_solid}
  Given $\Bar C$ as in Lemma~\ref{good_sextic} with contact conic $V$,
  there exists a quartic double solid such that $\Bar C$ arises as the
  discriminant locus of projection from one of the nodes and $V$ is
  the projection of the tangent cone. 
\end{proposition}
\begin{proof}
  We follow the construction on pages 104--105 of \cite{Izadi}. We take
  the double cover $\psi\colon Y\to\PP^2$ branched along $\Bar C$ and
  map it to $\PP^3$ by a linear system determined by $V$. 

  To define the linear system, take the desingularisation
  $\sigma\colon \Tilde Y\to Y$. The inverse image $\sigma^{-1}(V)$
  consists of two components, $V_+$ and $V_-$. We consider the linear
  systems $H\pm=|(\psi\sigma)^*\cO_{\PP^2}(1)\tensor \cO_{\Tilde
    Y}(V_\pm)|$. Either of these linear systems (and no others) maps
  the K3 surface $\Tilde Y$ onto a quartic surface $\Bar Y$, unique up
  to projective equivalence, with six nodes: five of these nodes are
  the images of the exceptional curves of $\sigma$, corresponding to
  nodes of $Y$, and the sixth is the image of $V_{\pm}$.  The discriminant
  curve under projection from this sixth node is $\Bar C$ up to
  projective equivalence. We refer to the proof of
  \cite[Theorem~2.1.1]{Izadi} for the computations of the degree and
  dimension of $H_\pm$ needed to justify these assertions. 

The quartic double solid whose existence is asserted in the theorem is
the double cover of $\PP^3$ branched along $\Bar Y$, and it is
immediate from the construction that it has the required properties. 
\end{proof}

\begin{corollary}\label{surjectivity_geom}
The map $\varrho\colon \cQ\to \cR_5$ is dominant. 
\end{corollary}
\begin{proof}
Immediate from Proposition~\ref{rebuild_double_solid}. 
\end{proof}
  
\section{Families of canonical curves}\label{cancurves}

In this section we give an alternative proof that $\theta$, and hence
$\varrho$, is dominant.  The method is to check directly, by
computation, that the Kodaira-Spencer map is locally surjective at a
test point. It does not rely on the special geometry of $\cQ$: it is a
method of checking computationally that a given family of $5$-nodal
sextics is general in the sense of moduli of genus $5$ curves. For
simplicity we assume in this section that $\KK$ is of characteristic zero. 

We first write down a local condition for $\theta$ to be dominant, {i.e.} 
generically surjective. 
\begin{lemma}\label{local_surjectivity}
  Let $u\colon X \to X'$ be a morphism, with $X'$
  irreducible. Then $u(X)$ is Zariski dense in
  $X'$ if and only if there exists a smooth point $P\in X$ such that the
  differential $d{u}_{P}\colon T_{X,P} \to T_{X',u(P)}$ is surjective. 
\end{lemma}

\begin{proof}
  Since $X'$ is irreducible the closure of $u(X)$ is
  $X'$ if and only if the dimension of one of
  its irreducible components is equal to $\dim X'$. Now it is
  enough to recall that the dimension of the irreducible component of
  $\overline{u(X)}$ containing a regular point $u(P)$ is given by the
  rank of the differential. 
\end{proof}
The tangent space to any scheme $X$ at a closed point
$P$ is the set of maps from $D =
\Spec{\KK[\varepsilon]/(\varepsilon^2)}$ to $X$ centred at $P$. 
For any morphism $u\colon X \to X'$ and any closed regular point
$P\in{X}$ the differential $d{u}_{P}\colon T_{X,P} \to T_{X',u(P)}$
is given by $\varphi
\mapsto u\circ\varphi$. 

Let $C$ be a canonically embedded curve of genus five, which we assume
to be given by the complete intersection of three quadrics in $\PP^4$
({i.e.}, by Petri's Theorem, non-trigonal). Two canonically embedded
curves of genus $g$ are isomorphic if and only if they are
projectively equivalent. 

We put $R_2=H^0(\cO_{\PP^4}(2))$, the degree~$2$ part of
$\KK[x_0,\dots,x_4]$, which we identify with the space of $4\times 4$ 
symmetric matrices over~$\KK$. 

The set of all canonical curves in $\PP^4$ is an open subset of the
Grassmannian $\Gr\big(3,R_2)$. Projective equivalence is then given by the
action of the group $\PGL(5)$ on $\PP^4$. 

But the Grassmannian itself can be realised as an orbit space, this
time under the action of $\GL(3)$, as follows. Let $V$ be the open set
inside the $45$-dimensional vector space ${R_2}\times{R_2}\times{R_2}$
where the three components span a $3$-dimensional subspace of $R_2$, 
and consider the action of $\GL(3)$ whose orbits are all the possible
bases for a given subspace. This is the action
\begin{equation*}
M(\qx{Q_1},\qx{Q_2},\qx{Q_3})_j =
 \sum_{i=1}^3m_{ji}\qx{Q_i},\qquad 1\le j\le 3,
\end{equation*}
where $M\in \GL(3)$ and $\phantom{}^t\underline{x}$ is the row vector
$(x_0,\dots,x_4)$, so $\qx{Q_i}\in R_2$ if $Q_i$ is a symmetric
matrix. 

The action of $N\in \PGL(5)$ is given by
\begin{equation*}
N(\qx{Q_1},\qx{Q_2},\qx{Q_3})_j =
 \qx{\phantom{,}^{t}\!NQ_jN},\qquad 1\le j\le 3. 
\end{equation*}
The two actions commute and we can regard one as acting on the orbit
space of the other. 

In order to investigate properties of a smooth family of deformations of a
canonical genus $5$ curve $C$ it is enough to consider the case in
which the base scheme is the spectrum of $A =
\KK[t_0,\ldots,t_n]/\gothm^2$, where $\gothm$ is the maximal ideal
generated by $t_0,\ldots,t_n$, corresponding (as a point of
$\Spec{A}$) to the curve $C$. Then the family is the scheme $\cC=
\Proj A[x_0,\ldots,x_4]/(F_1,F_2,F_3)$, where the coefficients of
$F_i$ depend linearly on the parameters $t_i$:
\begin{equation*}
F_i = H_i + \sum_{j=0}^nt_jH_{ij},
\end{equation*}
where $H_i$, $H_{ij}\in R_2$. The $n$ triples of quadrics
$(H_{1j},H_{2j},H_{3j})$ generate the linear subspace of
$\AA^{45}=R_2\times R_2\times R_2$ tangent to the family $\cC$. 

We want to compare this linear space with the tangent space to
$\cM_5$ at $C$. Our strategy is to work inside the
tangent space to $V$ at $s$: we construct a basis for all the trivial
deformations using the fact that they are
those given by the actions of $\PGL(5)$ and $\GL(3)$, and then check
how many of the above triples lie inside
this linear space.  

Around any point $v = (\qx{Q_1},\qx{Q_2},\qx{Q_3})$ in $V$ the
action of the two groups is linearised by the action of the
corresponding Lie algebras, so a system of generators for the linear
space tangent to the orbit passing through $v$ is simply determined by
applying a basis for the Lie algebra to it. The Lie algebra
$\mathfrak{gl}(3)$ is simply the whole space of three-by-three
matrices and its action is the same as the action of $\GL(3)$ so we
obtain a first set of trivial deformations given by the nine vectors
\begin{equation*}
(H_1,0,0), (H_2,0,0), \ldots, (0,0,H_2), (0,0,H_3). 
\end{equation*}

The algebra $\mathfrak{sl}(5)$ (which is the tangent space to $\PGL(5)$)
is the space of traceless $5 \times 5$ matrices, and its action is
determined as follows:
\begin{align*}
\phantom{}^t(N\underline{x})Q_i(N\underline{x})
        =& \phantom{}^t\underline{x}\phantom{\,}^t(I+\varepsilon\Delta)Q_i(I+\varepsilon\Delta)\underline{x} \\
        =& \phantom{}^t\underline{x}\big(Q_i + \varepsilon(\phantom{}^t\Delta{Q_i}+Q_i\Delta)\big)\underline{x} \\
        =& \qx{Q_i} + \varepsilon\underline{x}(\phantom{}^t\Delta{Q_i}+Q_i\Delta)\underline{x}. 
\end{align*}
Letting $\Delta$ vary among a basis for $\mathfrak{sl}(5)$ we get
another set of trivial deformations given by the $24$ vectors
\begin{equation*}
(\phantom{}^t\Delta{H_1}+Q_1\Delta,
 \phantom{}^t\Delta{H_2}+Q_2\Delta,
 \phantom{}^t\Delta{H_3}+Q_3\Delta   ). 
\end{equation*}

Now, given an $n$-dimensional family $\cF$ centred at $C$, we construct a matrix 
\begin{equation*}
M_{\cF}:= \left( \begin{array}{ccc}
 H_{11} & H_{21} & H_{31} \\
 H_{12} & H_{22} & H_{32} \\
 \vdots & \vdots & \vdots \\
 H_{1n} & H_{2n} & H_{3n} \\
 H_1 & 0 & 0 \\
 H_2 & 0 & 0 \\
 \vdots & \vdots & \vdots \\
 0 & 0 & H_2 \\
 0 & 0 & H_3 \\
 D_{21}{H_1}+H_1D_{12} &  D_{21}{H_2}+H_2D_{12} &  D_{21}{H_3}+H_3D_{12} \\
 D_{31}{H_1}+H_1D_{13} &  D_{31}{H_2}+H_2D_{13} &  D_{31}{H_3}+H_3D_{13} \\
 \vdots & \vdots & \vdots \\
 D_{55}{H_1}+H_1D_{55} &  D_{55}{H_2}+H_2D_{55} &  D_{55}{H_3}+H_3D_{55}
\end{array} \right)
\end{equation*}
The first $n$ rows are given by the family ${\cF}$: they are tangent
vectors at the central point $s = (H_1,H_2,H_3)$. The second set of
$9$ rows is given by the tangent vectors to the orbit of the
$\GL(3)$-action, and the last $24$ rows are the tangent vectors to the
orbits of the $\PGL(5)$-action described above. We have chosen a
vector space basis $D_{ij}$ for $\mathfrak{sl}(5)$, for example
$D_{ij}=\delta_{ij}$ for $i\neq j$ and $D_{ii}=\delta_{11}-\delta_{ii}$
for $1<i\le 4$. 

The linear space generated by the rows of $M_{\cF}$ is the
span inside the tangent space to $V$ of the three linear spaces
tangent respectively to the given family and to each of the two orbits
through $s$. To determine the dimension of this span we now need to
compute the rank of $M_{\cF}$. 

\begin{proposition}\label{surjectivity_criterion}
  Let $C$ be a smooth complete intersection of three linearly
  independent quadrics in $\PP^4$, and let $\cF$ be an $n$-dimensional
  family of deformations of $C$ as above. Suppose that $n\geq 12$. If
  the rank of $M_{\cF}$ is maximal then the Kodaira-Spencer map of
  $\cF$ at $C$ is surjective. 
\end{proposition}

\begin{proof}
  First observe that in the matrix $M_{\cF}$ there are $45$ columns,
  and under our assumptions there are at least $45$ rows. When the
  rank of the matrix $M_{\cF}$ is maximal the span of the three vector
  spaces, the two corresponding to trivial deformations and the one
  given by the family, is the whole of the tangent space to $V$ at the
  point $s$. Thus we are guaranteed the existence of enough linearly
  independent deformations, namely $12$, to fill the tangent space to
  $\cM_5$. 
\end{proof}

\begin{corollary}\label{surjectivity_comp}
The map $\varrho\colon \cQ\to \cR_5$ is dominant. 
\end{corollary}
\begin{proof}
This is now a straightforward computation of the rank of $M_\cF$ in
one particular case. We carried it out using {\it Macaulay\/}, with
points defined over a finite field (we chose $\FF_{101}$, for no
special reason). This is enough because if the rank is generically
maximal after reduction mod~$p$ it is also maximal in characteristic
zero. 

We chose the test point of $\cQ(\FF_{101})$ given by
\begin{align*}
u_2 &= 19x_0^2-33x_0x_1+50x_1^2-13x_0x_2+50x_1x_2-15x_2^2 \\
u_3 &= -2x_0^2x_1-35x_0x_1^2-18x_0^2x_2-8x_0x_1x_2-36x_1^2x_2-4x_0x_2^2+45x_1x_2^2 \\
u_4 &= -38x_0^2x_1^2-32x_0^2x_1x_2-32x_0x_1^2x_2-6x_0^2x_2^2-38x_0x_1x_2^2+2x_1^2x_2^2. 
\end{align*}
We arrived at this by first selecting six points
$P_0=(0\!:\! 0\!:\! 0\!:\! 1),\ldots,P_3=(1\!:\! 0\!:\! 0\!:\! 0), P_4=(1\!:\! 1\!:\! 1\!:\! 1),
P_5=(1\!:\! 2\!:\! 3\!:\! 4)\in\PP^3$, with ideals $\gothp_0,\ldots,\gothp_5$, to be
the prescribed nodes of a quartic and then choosing at random a
quartic $F=\in\gothp_0^2\cap\cdots\cap\gothp_5^2$. Then $u_2$, $u_3$ and
$u_4$ are defined by $F=u_2x_3^2+u_3x_3+u_4$, and the $6$-nodal
quartic surface is given by $F=0$. 

Having chosen $F$ at random one must check that it is suitably
general, namely that $X$ has no other singular points and that the
singularity of $X$ at each $P_i$ is a simple node. 

The $5$-nodal plane sextic in this example is given by $u_3^2=u_2u_4$,
with nodes at $(0\!:\! 0\!:\! 1)$, $(0\!:\! 1\!:\! 0)$, $(1\!:\! 0\!:\! 0)$, $(1\!:\! 1\!:\! 1)$ and
$(1\!:\! 2\!:\! 3)$. We construct the blowup of $\PP^2$ in these five points by
considering the linear system of cubics passing through them and from
this we can easily compute the canonical curve $\Tilde C$ as the
intersection of three quadrics $H_1$, $H_2$ and $H_3$. These at once
give us the last $33$ rows of $M_{\cQ}$. To compute the first $13$
rows one must know the family $\cQ$, that is $B_0$, near $X$. We can
obtain local coordinates on $B_0$ from the coordinates on
$\PP(I_4)\cross\PP^3$ by computing a Gr\"obner basis for the ideal of
$B_0$. The for each first-order deformation $X_j$ corresponding to a
local coordinate $t_j$ we compute the quadrics defining the canonical
curve $\Tilde C_j$ exactly as we did for $\Tilde C$. These quadrics
are $H_i+t_jH_{ij}$ (with a correct choice of coordinates) and we have
computed~$M_{\cQ}$. 
\end{proof}

\section{Conclusions}\label{conclusions}

We can now deduce our main result immediately from the results of
Sections~\ref{doublesolid} or \ref{cancurves}. 

\begin{theorem}\label{R5.unirat}
  The moduli space $\cR_5$ of \'etale double covers of curves
  of genus five is unirational. 
\end{theorem}

\begin{proof}
  This follows from Corollary~\ref{surjectivity_geom} or
  Corollary~\ref{surjectivity_comp}. 
\end{proof}

Theorem~\ref{R5.unirat} also provides a slightly different proof of a
theorem of Clemens~\cite{Clemens1}. 

\begin{corollary}\label{A4.unirat}
$\cA_4$ is unirational. 
\end{corollary}

This follows from Theorem~\ref{R5.unirat} because the Prym map
$p_5\colon \cR_5 \to \cA_4$ is dominant (see for instance
\cite{Beauville1}). The original proof of Clemens also starts from
quartic double solids, but Clemens exhibits the general principally
polarised abelian $4$-fold as an intermediate Jacobian rather than a
Prym variety.

\bigskip
\noindent
E.~Izadi\\
Department of Mathematics\\
University of Georgia\\
Athens\\ 
GA 30602-7403\\
USA\\ 
{\tt izadi@math.uga.edu}
\bigskip

\noindent
M.~Lo~Giudice,\ G.K.~Sankaran\\
Department of Mathematical Sciences\\
University of Bath\\
Bath BA2 7AY\\
England\\
{\tt marco.logiudice@gmail.com}\\
{\tt gks@maths.bath.ac.uk}

\end{document}